\documentclass[preprint,12pt]{elsarticle}

\usepackage{amssymb}
\usepackage{amsthm}
\usepackage{amsmath}
\usepackage{bm}
\usepackage{hyperref}

\numberwithin{equation}{section}
\newtheorem{thm}{Theorem}[section] 
\newtheorem{lem}{Lemma}[section]
\newtheorem{rmk}{Remark}[section]
\newtheorem{prop}{Proposition}[section]
\theoremstyle{definition}

\allowdisplaybreaks[4] 

\usepackage{lipsum}
\makeatletter
\def\ps@pprintTitle{%
 \let\@oddhead\@empty
 \let\@evenhead\@empty
 \def\@oddfoot{}%
 \let\@evenfoot\@oddfoot}
\makeatother

\begin{document}

\begin{frontmatter}



\title{Stability for a family of planar systems with nilpotent critical points}


\author[1]{Ziwei Zhuang \corref{cor1}}\ead{zhuangzw@mail2.sysu.edu.cn}
\author[1]{Changjian Liu}

\cortext[cor1]{Corresponding author}
\affiliation[1]{organization={School of Mathematics (Zhuhai), Sun Yat-sen University},
            city={Zhuhai},
            postcode={519086}, 
            country={P. R. China}}

\begin{abstract}
Consider a family of planar polynomial systems 
    $\dot x = y^{2l-1} - x^{2k+1},\\
    \dot y = -x +m y^{2s+1},$
    where
    $l,k,s\in\mathbb{N^*},$
    $2\le l \le 2s$
    and
    $m\in\mathbb{R}.$
    We study the center-focus problem on its origin which is a monodromic nilpotent critical point. 
    By directly calculating the generalized Lyapunov constants, 
    we find that the origin is always a focus 
    and we complete the classification of its stability.
    This includes the most difficult case: $s=kl$ and $m=(2k+1)!!/(2kl+1)!_{(2l)}.$
    In this case, we prove that the origin is always unstable. 
    Our result extends and completes a previous one.
\end{abstract}


\end{frontmatter}


\section{Introduction}

One of the classic and important researches on qualitative properties of planar analytic systems
is the center-focus problem on a monodromic critical point (the origin in general).
It is known that if the origin is monodromic and has a non-zero Jacobian matrix, 
the eigenvalues of the matrix are either
a pair of pure imaginary numbers,
i.e., a linear center or focus, 
or two zeros, i.e., a nilpotent critical point.
For the linear case,
the center-focus problem can be solve by algorithms 
that calculate the Lyapunov constants in complex or polar coordinates,
see for instance \cite[Chapter 4]{dumortier2006qualitative}.
For the nilpotent case, the monodromy further requires that 
some parameters of the system lie in a suitable region.
Details can be referred to \cite{andreev1958investigation}
or \cite[Chapter 3]{dumortier2006qualitative}.
The center condition for a nilpotent origin is more complicated.
Motivated by the linear case,  there is also a procedure 
to verify the nilpotent center by calculating
the generalized Lyapunov constants defined by the generalized trigonometric functions,
see \cite{alvarez2005monodromy,lyapunov1966stability,gasull2001new}.
Another popular approach is to change the original system to the normal form
(the form of the Li{\'e}nard system) and to check whether some function is odd.
Details can be referred to \cite{moussu1982symetrie}, 
and also \cite{chen2022new} for some extensions.
It is worth mentioning that once the monodromic origin is verified 
not to be a center, the following is to consider its stability.
However, in some cases these two problems can be solved simultaneously.
For example, the monodromic origin is a focus if and only if 
there is a non-vanishing Lyapunov constant, and
the sign of the first non-vanishing Lyapunov constant
determines the stability of the focus.

This paper deals with a certain family of nilpotent systems 
which is motivated from the following ones
\begin{align}\label{generailzed system}
    \left\{
    \begin{aligned}
    &\dot x = y^3 - x^{2k+1},\\
    &\dot y = -x +m y^{2s+1}, 
    \quad m\in\mathbb R,\; s,k\in\mathbb{N^*}.
    \end{aligned}
    \right.
\end{align}
It is known that the origin
is always a  monodromic nilpotent critical point.
When 
$s\neq 2k,$ 
Galeotti and Gori \cite{galeotti1989bifurcations} showed that the origin is a focus
and they clarified its stability.
For the rest case,  
García-Saldaña et al. \cite{garcia2014bifurcation}
went further and showed that the stability changes at
\begin{equation*}
    m=\frac{(2k+1)!!}{(4k+1)!_{(4)}}.
\end{equation*}
Hereby the notation $a!_{(n)}$ 
with non-negative real number $a$ and positive integer $n$
are recursively defined by
\begin{align*}
    a!_{(n)}=\left\{
    \begin{aligned}
        &a\cdot(a-n)!_{(n)},\quad &a>0, 
        \\
        &1,\quad &a\le0,
    \end{aligned}
    \right.
\end{align*}
and the cases $n=1$ and $2$ are abbreviated as $a!$ and $a!!,$ respectively.
Precisely, the authors in \cite{garcia2014bifurcation} 
used the method that calculates the generalized Lyapunove constants
and they gave the following result.
\begin{thm}{\cite[Theorem 1.3]{garcia2014bifurcation}}
\label{Theorem: garcia's result}
    Consider system \eqref{generailzed system}.
    \begin{enumerate}[(i)]
        \item When $s<2k,$ the origin is an attractor for
        $m\le0$ and a repeller for $m>0.$
        \item When $s>2k,$ the origin is always an attractor.
        \item When $s=2k,$ the origin is an attractor for
        $m<(2k+1)!!/(4k+1)!_{(4)}$ 
        and a repeller for $m>(2k+1)!!/(4k+1)!_{(4)}.$ 
        Moreover, when $k=1$ and $m=3/5$ the origin is a repeller.
    \end{enumerate}
\end{thm}
Recently Chen et al. \cite{chen2022new}
introduced a new sufficient condition to determine the stability of the nilpotent focus
based on the intermediate system 
which is transformed between the Jordan normal form and the Li\'enard normal form.
As an application on system \eqref{generailzed system},  
they also provided a result similar to Theorem \ref{Theorem: garcia's result},
with less calculations.

Apparently, Theorem \ref{Theorem: garcia's result} does not cover
all the cases. When $s=2k>2$ and $m=m^*,$ 
the behavior of the origin is still indeterminate. 
For this case, the center-focue problem was sovled by Caubergh \cite[Theorem 25]{caubergh2015bifurcation} in 2015. 
The author found that the assumption of a period annulus of the origin contradicts to its outer bound being a hyperbolic 2-saddle cycle.
Thus, the origin must be a focus.

In this paper, we consider a more general family of systems
\begin{equation}\label{main system}
    \left\{
  \begin{aligned}
  &\dot x = y^{2l-1} - x^{2k+1},\\
  &\dot y = -x +m y^{2s+1}, \quad l,k,s\in\mathbb{N^*},\;
  2\le l \le 2s.
  \end{aligned}
  \right.
\end{equation}
The restricted range of $l$ ensures that the origin is always a monodromic nilpotent critical point.
We only concentrate on the center-focus problem on the origin
and determine the stability if it is a focus.
We follow the classic method presented in \cite{garcia2014bifurcation},
that is to calculate the first non-vanishing generalized Lyapunov constant.
Details will be introduced in the next section.
The first result is a simple extension of Theorem \ref{Theorem: garcia's result}.
Define
$$ m^*\triangleq\frac{(2k+1)!!}{(2kl+1)!_{(2l)}}. $$
\begin{thm}\label{Theorem_extension for previous one}
     Consider system \eqref{main system}.
    \begin{enumerate}[(i)]
        \item When $s<kl,$ the origin is an attractor for
        $m\le0$ and a repeller for $m>0.$
        \item When $s>kl,$ the origin is always an attractor.
        \item When $s=kl,$ the origin is an attractor for
        $m<m^*$ and a repeller for $m>m^*.$ 
    \end{enumerate}
\end{thm}
What really interests us is the rest case:
$s=kl$ and $m=m^*.$
By the same method with a mass of tedious calculations, 
we obtain our main result.
\begin{thm}\label{Theorem_main result: unstable focus}
When $s=kl$ and $m=m^*,$
the origin of system \eqref{main system} is always a repeller.
\end{thm}
Obviously, when $l=2,$ 
Theorem \ref{Theorem_main result: unstable focus} completes Theorem \ref{Theorem: garcia's result}.

\section{The first non-vanishing generalized Lyapunov constant
and the proof of Theorem \ref{Theorem_extension for previous one}}

We first introduce the generalized polar coordinates
from the solution of the initial value problem
\begin{equation*}
\begin{aligned}
  &\frac{\mathrm{d} x}{\mathrm d \theta}=-y^{2p-1},
  \quad
  \frac{\mathrm{d} y}{\mathrm d \theta}=x^{2q-1},
  \\
  &x(0)=\sqrt[2q]{1/p},
  \quad
  y(0)=0.
\end{aligned}
\end{equation*}
Its solution, denoted by
$$
x(\theta)=\text{Cs}(\theta),
\quad
y(\theta)=\text{Sn}(\theta),
$$ are called the
generalized trigonometric functions
which are introduced by Lyapunov in \cite{lyapunov1966stability}.
They are period functions with period
$$
\Omega=2p^{-\frac{1}{2q}}q^{-\frac{1}{2p}}
\frac{\Gamma\left(\frac{1}{2p}\right)
  \Gamma\left(\frac{1}{2q}\right)}
  {\Gamma\left(\frac{1}{2p}+\frac{1}{2q}\right)},
$$
and they satisfy
$p\text{Cs}^{2q}(\theta)+q\text{Sn}^{2p}(\theta)=1.$
\textbf{In this paper, we always take}
$\bm{p=1}$
\textbf{and}
$\bm{q=l}.$
The next result is from \cite[Lemma 5]{gasull2001new}.
\begin{lem}\label{Lemma: integrals of Sn and Cs}
Let
$i$
and
$j$
be fixed natural numbers. The following items hold:
\begin{enumerate}[(i)]
  \item
  $
  \int\textnormal{Sn}\theta\textnormal{Cs}^j\theta\mathrm{d}\theta
  =-\frac{\textnormal{Cs}^{j+1}}{j+1} + \textnormal{Const.};
  $
  \item
  $
  \int\textnormal{Sn}^i\theta\textnormal{Cs}^{2l-1}\theta\mathrm{d}\theta
  =\frac{\textnormal{Sn}^{i+1}}{i+1} + \textnormal{Const.};
  $
  \item
  $
  \int\textnormal{Sn}^i\theta\textnormal{Cs}^j\theta\mathrm{d}\theta
  =-\frac{\textnormal{Sn}^{i-1}\theta\textnormal{Cs}^{j+1}\theta}{(i-1)l+j+1}
  +\frac{i-1}{(i-1)l+j+1}
  \int\textnormal{Sn}^{i-2}\theta\textnormal{Cs}^j\theta\mathrm{d}\theta;
  $
  \item
  $
  \int\textnormal{Sn}^i\theta\textnormal{Cs}^j\theta\mathrm{d}\theta
  =\frac{l\textnormal{Sn}^{i+1}\theta\textnormal{Cs}^{j-2l+1}\theta}{(i-1)l+j+1}
  +\frac{j-2l+1}{(i-1)l+j+1}
  \int\textnormal{Sn}^{i}\theta\textnormal{Cs}^{j-2l}\theta\mathrm{d}\theta;
  $
  \item
  $
  \int_0^{\Omega}
  \textnormal{Sn}^i\theta\textnormal{Cs}^j\theta\mathrm{d}\theta=0,
  $
  when either
  $i$
  or
  $j$
  is odd;
  \item
  $
  \int_0^{\Omega}
  \textnormal{Sn}^i\theta\textnormal{Cs}^j\theta\mathrm{d}\theta
  =\frac{2}{l^{\frac{i+1}{2}}}
  \frac{\Gamma\left(\frac{i+1}{2}\right)
  \Gamma\left(\frac{j+1}{2l}\right)}
  {\Gamma\left(\frac{i+1}{2}+\frac{j+1}{2l}\right)},
  $
  when both
  $i$
  and
  $j$
  are even.
\end{enumerate}
\end{lem}

By applying the transformation 
$(x,y)\rightarrow
\left(r^l\text{Sn}(\theta),r\text{Cs}(\theta)\right)$
to system \eqref{main system}, it follows that
\begin{equation*}
  \left\{
  \begin{aligned}
  &\dot r = mr^{2s+1}\text{Cs}^{2s+2l}(\theta)
        -r^{2kl+1}\text{Sn}^{2k+2}(\theta),\\
  &\dot \theta = r^{l-1} - lm r^{2s}\text{Sn}(\theta)\text{Cs}^{2s+1}(\theta)
        -r^{2kl}\text{Sn}^{2k+1}(\theta)\text{Cs}(\theta) .
  \end{aligned}
  \right.
\end{equation*}
In a small neighbourhood of the origin, the orbits are characterized by
the equation
\begin{equation}\label{mainsystem ode}
    \frac{\mathrm{d}r}{\mathrm{d}\theta}
    =\frac{ mr^{2s-l+2}\text{Cs}^{2s+2l}(\theta)-r^{(2k-1)l+2}\text{Sn}^{2k+2}(\theta) }
    { 1 -lm r^{2s-l+1}\text{Sn}(\theta)\text{Cs}^{2s+1}(\theta)
    - r^{(2k-1)l+1}\text{Sn}^{2k+1}(\theta)\text{Cs}(\theta) }.
\end{equation}
Denote its solution by 
$r=r(\theta,\rho)$
with
$r(\theta,0)\equiv0$
and
$r(0,\rho)=\rho$
for sufficiently small 
$\rho\ge0.$
Assuming $r(\theta,\rho)$ in the following form
\begin{equation}\label{r(theta,rho) in power series}
  r(\theta,\rho)=\rho+
    \sum_{i=2}^{+\infty}u_{i}(\theta)\rho^i
  \quad\text{with}\quad u_i(0)=0,
\end{equation}
the generalized Lyapunov constant is given by
$u_i(\Omega).$ 
Apparently, the sign of the first non-vanishing 
$u_i(\Omega)$  
decides the stability of the origin of system \eqref{main system}.
To calculate 
$u_i(\theta),$
we expand the right-hand side of \eqref{mainsystem ode} 
in the power series at $r=0.$
When $s\neq kl,$ we can see that
\begin{equation}\label{mainsystem ode expansion: s≠kl}
  \frac{\mathrm{d}r}{\mathrm{d}\theta}
  =\left\{
  \begin{aligned}
  & -r^{(2k-1)l+2}\text{Sn}^{2k+2}(\theta) + O\left(r^{(2k-1)l+3}\right),
  \quad && s>kl \quad\text{or}\quad m=0,
  \\
  & m r^{2s-l+2}\text{Cs}^{2s+2l}(\theta) + O\left(r^{2s-l+3}\right),
  \quad && s<kl \quad \text{and}\quad m\neq0.
  \end{aligned}
  \right.
\end{equation}
Substituting
\eqref{r(theta,rho) in power series}
into \eqref{mainsystem ode expansion: s≠kl}
and equating the coefficient of each $\rho^i$ at both sides,
the first non-vanishing $u_i(\theta)$ is
\begin{align*}
    &u_{(2k-1)l+2}(\theta)=-\int_{0}^{\theta}\text{Sn}^{2k+2}(\xi)\mathrm{d}\xi
    \quad \text{if}\quad s>kl\quad\text{or}\quad m=0,
    \\
    \text{or}\quad
    &u_{2s-l+2}(\theta)=m\int_{0}^{\theta}\text{Cs}^{2s+2l}(\xi)\mathrm{d}\xi
    \quad \text{if}\quad s<kl \quad\text{and}\quad m\neq0.
\end{align*}
Obviously, it follows that
$u_{(2k-1)l+2}(\Omega)<0$
if $s>kl$ or $m=0,$
$u_{2s-l+2}(\Omega)<0$ if $s<kl$ and $m<0,$
and
$u_{2s-l+2}(\Omega)>0$ if $s<kl$ and $m>0.$
These imply Theorem \ref{Theorem_extension for previous one} (i) and (ii).

When $s=kl,$ Eq. \eqref{mainsystem ode} becomes
\begin{equation*}
  \frac{\mathrm{d}r}{\mathrm{d}\theta}
  = \frac{ v(\theta) r^{(2k-1)l+2} }
   {1 - w(\theta) r^{(2k-1)l+1} },
\end{equation*}
where
\begin{equation*}
\begin{aligned}
v(\theta)=&m\text{Cs}^{2(k+1)l}(\theta)-\text{Sn}^{2k+2}(\theta),
\\
w(\theta)=&lm\text{Sn}(\theta)\text{Cs}^{2kl+1}(\theta)+
        \text{Sn}^{2k+1}(\theta)\text{Cs}(\theta) .
\end{aligned}
\end{equation*}
Its expansion at $r=0$ is 
\begin{equation}\label{mainsystem ode expansion: s=kl}
  \frac{\mathrm{d}r}{\mathrm{d}\theta}
  =v(\theta)r^{K+1} + v(\theta)w(\theta)r^{2K+1}
   +v(\theta)w^2(\theta)r^{3K+1}+O\left(r^{4K+1}\right),
\end{equation}
where $K=(2k-1)l+1.$
\begin{lem}\label{Lemma ui=0}
Let $n$ be any natural number.
If
$
j\neq nK+1,
$
then
$u_j(\theta)\equiv0.$
\begin{proof}
Substituting
\eqref{r(theta,rho) in power series}
into \eqref{mainsystem ode expansion: s=kl},
it is easy to see that the least order of
$\rho$
in the right-hand side of \eqref{mainsystem ode expansion: s=kl} is
$K+1,$
which implies that
$u_j(\theta)\equiv0$
for 
$j=2,\ldots,K.$

Inductively, we assume that
$u_j(\theta)\equiv0$
for
$
j\in\left\{1,2,\ldots,nK+1\right\}
$ \\ $
/\left\{1,K+1,\ldots,iK+1,
\ldots,nK+1\right\}
$
with some natural number
$n.$
It is to say that \eqref{r(theta,rho) in power series}
becomes 
\begin{equation*}
  r(\theta,\rho)=\rho + u_{K+1}(\theta)\rho^{K+1}
  +u_{2K+1}(\theta)\rho^{2K+1}
  +\cdots
  +u_{nK+1}(\theta)\rho^{nK+1}
  +O\left(\rho^{nK+2}\right).
\end{equation*}
Substituting it into the right-hand side of
\eqref{mainsystem ode expansion: s=kl}, we should focus on
whether the monomial
$\rho^{nK+h+1}$
vanishes for 
$h=1,\ldots,K-1.$
If any 
$\rho^{nK+h+1}$
is contributed by 
$r^{jK+1}$
for some integer
$j \ge 1,$
then there must be 
$n+h+1$
nonnegative integers
$\alpha_0,\alpha_1,\ldots,\alpha_n,
\beta_1,\ldots,\beta_h$
such that
\begin{align*}
  &\alpha_0+\alpha_1+\cdots+\alpha_n
  +\beta_1+\cdots+\beta_h=jK+1,
  \\
  &\alpha_0+(K+1)\alpha_1+\cdots+(nK+1)\alpha_n
  +(nK+2)\beta_1+
  \cdots+(nK+h+1)\beta_h
  \\
  &=nK+h+1.
\end{align*}
The second equation minus the first one becomes
\begin{equation*}
  K\left( \sum_{i=1}^{n}i\alpha_i+n\sum_{i=1}^{h}\beta_i \right)
  +\sum_{i=1}^{h}i\beta_i=(n-j)K+h,
\end{equation*}
which implies that $K \mid h-\sum_{i=1}^{h}i\beta_i.$
Since 
$0<h<K,$
it follows that not all $\beta_i$'s vanish and
$h-\sum_{i=1}^{h}i\beta_i \le 0.$
Then, we have
\begin{equation*}
  nK>(n-j)K
  =K\left( \sum_{i=1}^{n}i\alpha_i+n\sum_{i=1}^{h}\beta_i \right)
  +\sum_{i=1}^{h}i\beta_i-h
  \ge nK,
\end{equation*}
which is a contradiction.

Consequently, none of 
$r^{jK+1}$
can contribute a monomial 
$\rho^{nK+h+1}.$
That is to say
$u_{nK+h+1}(\theta)\equiv0$
for 
$h=1,\ldots,K-1,$
which completes the proof by induction.
\end{proof}
\end{lem}

By Lemma \ref{Lemma ui=0}, we then calculate
$u_{nK+1}(\theta)$ 
to obtain the first non-vanishing Lyapunov constant.
It is to see that
\begin{equation*}
    u_{K+1}(\theta)=\int_0^\theta v(\xi)\mathrm{d}\xi
    =m\int_0^{\theta}\text{Cs}^{2(k+1)l}(\xi)\mathrm{d}\xi
     -\int_0^{\theta}\text{Sn}^{2k+2}(\xi)\mathrm{d}\xi.
\end{equation*}
By Lemma \ref{Lemma: integrals of Sn and Cs} (vi),
it follows that the sign of $u_{K+1}(\Omega)$ changes at 
$$
m=\frac{\int_0^{\Omega}\text{Sn}^{2k+2}(\theta)\mathrm{d}\theta}
{\int_0^{\Omega}\text{Cs}^{2(k+1)l}(\theta)\mathrm{d}\theta}
=\frac{ \frac{2}{l^{ k+1+\frac{1}{2} }} \frac{\Gamma\left(k+1+\frac{1}{2}\right)
        \Gamma\left(\frac{1}{2l}\right)}
        {\Gamma\left(k+1+\frac{1}{2}+\frac{1}{2l}\right)}
        }{
        \frac{2}{l^{ \frac{1}{2} }} 
        \frac{\Gamma\left( \frac{1}{2} \right)
        \Gamma\left(k+1+\frac{1}{2l}\right)}
        {\Gamma\left(k+1+\frac{1}{2}+\frac{1}{2l}\right)}
        }
=\frac{1}{l^{k+1}}
\frac{\left(k+\frac{1}{2}\right)!}{\left(k+\frac{1}{2l}\right)!}
=m^*.
$$
Precisely, $u_{K+1}(\Omega)<0$ when $m<m^*,$
and $u_{K+1}(\Omega)>0$ when $m>m^*.$
These imply Theorem \ref{Theorem_extension for previous one} (iii)
and the proof of Theorem \ref{Theorem_extension for previous one} is completed.

Finally, when $s=kl$ and $m=m^*,$
we have to calculate more $u_{nK+1}(\theta)$ 
to obtain the first non-vanishing generalized Lyapunov constant.
Fortunately the calculation stops at 
$n=3,$
and we obtain
\begin{align*}
    u_{2K+1}(\theta)=&\frac{1}{2}(K+1)
    \left( \int_0^\theta v(\xi)\mathrm{d}\xi \right)^2
    + \int_0^\theta v(\xi)w(\xi)\mathrm{d}\xi,
    \\
    u_{3K+1}(\theta)=&\frac{1}{6}(K+1)(2K+1)
    \left( \int_0^\theta v(\xi)\mathrm{d}\xi \right)^3
    \\
    &- K\int_0^\theta v(\xi) 
    \left( \int_0^\xi v(\eta)w(\eta)\mathrm{d}\eta \right)
    \mathrm{d}\xi
    + \int_0^\theta v(\xi)w^2(\xi)\mathrm{d}\xi.
\end{align*}
Recall that
$\int_0^\Omega v(\theta)\mathrm{d}\theta=0$
when $m=m^*,$
and by Lemma \ref{Lemma: integrals of Sn and Cs} (v),
$\int_0^\Omega v(\theta)w(\theta)\mathrm{d}\theta=0.$
It follows that the first possibly non-vanishing generalized Lyapunov constant is
$$
u_{3K+1}(\Omega)
=-K\int_0^\Omega v(\theta) 
\left(\int_0^\theta v(\xi)w(\xi)\mathrm{d}\xi\right)
\mathrm{d}\theta
+ \int_0^\Omega v(\theta)w^2(\theta)\mathrm{d}\theta.
$$
Recall that $K=(2k-1)l+1.$
The proof of Theorem \ref{Theorem_main result: unstable focus}
is to verify $u_{3K+1}(\Omega)>0$ for all $l$ and $k.$

\section{Proof of Theorem \ref{Theorem_main result: unstable focus}}

Let 
\begin{equation}
  V=-\int_0^{\Omega}v(\theta)
  \left(\int_0^{\theta}v(\xi)w(\xi)\mathrm{d}\xi\right)
  \mathrm{d}\theta
  \quad\text{and}\quad
  W=\int_0^{\Omega}v(\theta)w^2(\theta)\mathrm{d}\theta.
\end{equation}
We are going to show 
$V>0$
and
$W>0,$
respectively.
Before this, we give a simple result about 
a kind of infinite products that
will be frequently used in this section.

One knows that for a real series 
$\{a_n\}_{n=1}^{+\infty}$
with $a_n\ge0,$
the infinite product 
$\prod_{n=1}^{+\infty}\left(1+a_n\right)$
converges if and only if
$\sum_{n=1}^{+\infty}a_n$
converges. Similar statement holds for
$-1<a_n\le0.$
In particular, when 
$a_n=O(1/n^2),$ 
we have the following estimates.

\begin{lem}\label{Lemma: infinite product converge}
Let 
$\{a_n\}_{n=1}^{+\infty}$
be a real series such that 
$0\le a_n<1$
and
$\lim_{n\rightarrow+\infty} n^2 a_n=\bar a<+\infty.$
Then for sufficiently small
$\varepsilon>0,$
there is a positive integer $N_{\varepsilon}$
such that
\begin{equation}\label{general infinite product 1}
  \prod_{n=1}^{+\infty}\left(1+a_n\right)
  <
  \exp\left(\frac{\bar a + \varepsilon}{N_{\varepsilon}}\right)
\prod_{n=1}^{N_{\varepsilon}}\left(1+a_n\right),
\end{equation}
or
\begin{equation}\label{general infinite product 2}
  \prod_{n=1}^{+\infty}\left(1-a_n\right)
  >
  \exp\left(-\frac{\bar a + \varepsilon}{N_{\varepsilon}}\right)
\prod_{n=1}^{N_{\varepsilon}}\left(1-a_n\right).
\end{equation}

\begin{proof}
For any $\varepsilon>0,$ take
$N_{\varepsilon}$ to satisfy
$$a_n<\frac{\bar a +\varepsilon}{n^2},\quad
\forall\;n>N_{\varepsilon}.$$
Then, we have for any $n>N_\varepsilon$
\begin{align*}
    \prod_{i={N_{\varepsilon}+1}}^{n}(1+a_i)
<&\exp\left(\sum_{i=N_{\varepsilon}+1}^{n}a_i\right)
<\exp\left(\sum_{i=N_{\varepsilon}+1}^{n}
\frac{\bar a +\varepsilon}{i^2}\right)
\\
<&\exp\left(\sum_{i=N_{\varepsilon}+1}^{n}
\frac{\bar a +\varepsilon}{(i-1)i}\right)
<\exp\left(\frac{\bar a +\varepsilon}{N_{\varepsilon}}\right).
\end{align*}
This gives $\eqref{general infinite product 1},$
which holds for $a_n\ge0.$
On the other hand, since $a_n$ is restricted on
$[0,1)$, notice that
$$\frac{1}{1-a_n}=1+\frac{a_n}{1-a_n},
\quad 
\frac{a_n}{1-a_n}\ge0
\quad\text{and}\quad
\frac{n^2a_n}{1-a_n}\rightarrow\bar a.$$
Then \eqref{general infinite product 2} is obtained by
applying \eqref{general infinite product 1} on 
$1/\prod_{n=1}^{+\infty}(1-a_n).$
In this case, $N_\varepsilon$ is taken to satisfy
$$
\frac{a_n}{1-a_n}
<\frac{\bar a +\varepsilon}{n^2},\quad
\forall\;n>N_{\varepsilon}.
$$
\end{proof}
\end{lem}

\subsection{Proof of $V>0$}

By direct calculation, we divide $V$ into four parts:
\begin{equation*}
    V=V_{1}+V_{2}+V_{3}+V_{4},
\end{equation*}
where
\begin{align*}
  V_{1}=&-l(m^*)^2\int_0^{\Omega} v(\theta)
  \left(
  \int_0^{\theta}\text{Sn}(\xi)\text{Cs}^{2(2k+1)l+1}(\xi)\mathrm{d}\xi
  \right)\mathrm{d}\theta,
  \\
  V_{2}=&\int_0^{\Omega} v(\theta)
  \left(
  \int_0^{\theta}
  \text{Sn}^{4k+3}(\xi)\text{Cs}(\xi)\mathrm{d}\xi
  \right)\mathrm{d}\theta,
  \\
  V_{3}=&-m^*\int_0^{\Omega}v(\theta)
  \left(
  \int_0^{\theta}
  \text{Sn}^{2k+1}(\xi)\text{Cs}^{2l(k+1)+1}(\xi)\mathrm{d}\xi
  \right)\mathrm{d}\theta,
  \\
  V_{4}=&lm^*\int_0^{\Omega}v(\theta)
  \left(
  \int_0^{\theta}
  \text{Sn}^{2k+3}(\xi)\text{Cs}^{2kl+1}(\xi)\mathrm{d}\xi
  \right)\mathrm{d}\theta.
\end{align*}
By recursively using
Lemma \ref{Lemma: integrals of Sn and Cs} (iii),
we have
\begin{align*}
  \int 
  \text{Sn}^{4k+3}(\theta)\text{Cs}(\theta)
  \mathrm{d}\theta
  =&-\frac{\text{Cs}^2(\theta)}{2(2k+1)l+2}
  \frac{(2k+1)!}{(2kl+1)!_{(l)}}
  \\
  &\times
  \sum_{i=0}^{2k+1}
  \frac{((i-1)l+1)!_{(l)}}{i!}
  \text{Sn}^{2i}(\theta)
  +\text{Const.,}
  \\
  \int
  \text{Sn}^{2k+1}(\theta)\text{Cs}^{2(k+1)l+1}(\theta)
  \mathrm{d}\theta
  =&-\frac{\text{Cs}^{2(k+1)l+2}(\theta)}{2(2k+1)l+2}
  \frac{k!}{(2kl+1)!_{(l)}}
  \\
  &\times
  \sum_{i=0}^{k}
  \frac{((k+i)l+1)!_{(l)}}{i!}
  \text{Sn}^{2i}(\theta)
  +\text{Const.},
  \\
  \int
  \text{Sn}^{2k+1}(\theta)\text{Cs}^{2(k+1)l+1}(\theta)
  \mathrm{d}\theta
  =&-\frac{\text{Cs}^{2kl+2}(\theta)}{2(2k+1)l+2}
  \frac{(k+1)!}{(2kl+1)!_{(l)}}
  \\
  &\times
  \sum_{i=0}^{k+1}
  \frac{((k+i-1)l+1)!_{(l)}}{i!}
  \text{Sn}^{2i}(\theta)
  +\text{Const.}
\end{align*}
Recall that
$\int_0^{\Omega}v(\theta)\mathrm{d}\theta=0.$
Then,
$V_{j}$ can be denoted by
\begin{equation*}
  V_{j}=(-1)^{j+1}\sum_i \left(b_{j}^i - c_{j}^i\right),
\end{equation*}
where
\begin{align*}
  b_1^0=&\frac{l(m^*)^3}{2(2k+1)l+2}\int_0^{\Omega}
  \text{Cs}^{2(3k+2)l+2}(\theta)\mathrm{d}\theta,
  \\
  c_1^0=&\frac{l(m^*)^2}{2(2k+1)l+2}\int_0^{\Omega}
  \text{Sn}^{2(k+1)}
  \text{Cs}^{2(2k+1)l+2}(\theta)\mathrm{d}\theta,
  \\
  b^i_{2}=&\frac{m^*}{2(2k+1)l+2}
   \frac{(2k+1)!}{(2kl+1)!_{(l)}}
  \frac{((i-1)l+1)!_{(l)}}{i!}
  \int_0^{\Omega}
  \text{Sn}^{2i}(\theta)
  \text{Cs}^{2(k+1)l+2}(\theta)
  \mathrm{d}\theta,
  \\
  c^i_{2}=&\frac{1}{2(2k+1)l+2}
   \frac{(2k+1)!}{(2kl+1)!_{(l)}}
  \frac{((i-1)l+1)!_{(l)}}{i!}
  \int_0^{\Omega}
  \text{Sn}^{2(k+i+1)}(\theta)
  \text{Cs}^{2}(\theta)
  \mathrm{d}\theta,
  \\
  b^i_{3}=&\frac{(m^*)^2}{2(2k+1)l+2}
   \frac{k!}{(2kl+1)!_{(l)}}
  \frac{((k+i)l+1)!_{(l)}}{i!}
  \int_0^{\Omega}
  \text{Sn}^{2i}(\theta)
  \text{Cs}^{4(k+1)l+2}(\theta)
  \mathrm{d}\theta,
  \\
  c^i_{3}=&\frac{m^*}{2(2k+1)l+2}
   \frac{k!}{(2kl+1)!_{(l)}}
  \frac{((k+i)l+1)!_{(l)}}{i!}
  \int_0^{\Omega}
  \text{Sn}^{2(k+i+1)}(\theta)
  \text{Cs}^{2(k+1)l+2}(\theta)
  \mathrm{d}\theta,
  \\
  b^i_{4}=&\frac{l(m^*)^2}{2(2k+1)l+2}
   \frac{(k+1)!}{(2kl+1)!_{(l)}}
  \frac{((k+i-1)l+1)!_{(l)}}{i!}
  \int_0^{\Omega}
  \text{Sn}^{2i}(\theta)
  \text{Cs}^{2(2k+1)l+2}(\theta)
  \mathrm{d}\theta,
  \\
  c^i_{4}=&\frac{lm^*}{2(2k+1)l+2}
   \frac{(k+1)!}{(2kl+1)!_{(l)}}
  \frac{((k+i-1)l+1)!_{(l)}}{i!}
  \int_0^{\Omega}
  \text{Sn}^{2(k+i+1)}(\theta)
  \text{Cs}^{2kl+2}(\theta)
  \mathrm{d}\theta.
\end{align*}
Now, $V$ becomes
\begin{equation*}
    V=b_1^0-c_1^0
     -\sum_{i=0}^{2k+1}\left(b_{2}^i-c_{2}^i\right)
    + \sum_{i=0}^{k}\left(b_{3}^i-c_{3}^i\right)
    - \sum_{i=0}^{k+1}\left(b_{4}^i-c_{4}^i\right).
\end{equation*}
Each of its terms can be solved by Lemma \ref{Lemma: integrals of Sn and Cs} (vi).

\begin{prop} 
\label{Lemma: estimate of infinitesimal}
\begin{enumerate}[(i)]
    \item $b_3^i-b_4^{i+1}>0;$
    \item $-c_3^{i}+c_4^{i+1}>0;$
    \item $-c_1^0-\sum_{i=1}^{2k+1}b_2^i +c_2^0>0.$
\end{enumerate}
\begin{proof}
    (i) For $i\ge0,$ we have
    \begin{align*}
        \frac{b_4^{i+1}}{b_3^{i}}=& l\frac{k+1}{i+1}
        \frac{ \frac{2}{l^{ i+1+\frac{1}{2} }} \frac{\Gamma\left(i+1+\frac{1}{2}\right)
        \Gamma\left(2k+1+\frac{3}{2l}\right)}
        {\Gamma\left(2k+i+2+\frac{1}{2}+\frac{3}{2l}\right)}
        }{
        \frac{2}{l^{ i+\frac{1}{2} }} 
        \frac{\Gamma\left( i+\frac{1}{2} \right)
        \Gamma\left(2k+2+\frac{3}{2l}\right)}
        {\Gamma\left(2k+i+2+\frac{1}{2}+\frac{3}{2l}\right)}
        }
        =\frac{\left(i+\frac{1}{2}\right)(k+1)}
        {(i+1)(2k+1+\frac{3}{2l})}<1.
    \end{align*}
    (ii) Recall that $l\ge2.$ For $i\ge0,$ we have
    \begin{align*}
        \frac{c_4^{i+1}}{c_3^{i}}=& l\frac{k+1}{i+1}
        \frac{ \frac{2}{l^{ k+i+2+\frac{1}{2} }} \frac{\Gamma\left(k+i+2+\frac{1}{2}\right)
        \Gamma\left(k+\frac{3}{2l}\right)}
        {\Gamma\left(2k+i+2+\frac{1}{2}+\frac{3}{2l}\right)}
        }{
        \frac{2}{l^{ k+i+1+\frac{1}{2} }} 
        \frac{\Gamma\left( k+i+1+\frac{1}{2} \right)
        \Gamma\left(k+1+\frac{3}{2l}\right)}
        {\Gamma\left(2k+i+2+\frac{1}{2}+\frac{3}{2l}\right)}
        }
        =\frac{\left(k+i+\frac{3}{2}\right)(k+1)}
        {(i+1)\left(k+\frac{3}{2l}\right)}>1.
    \end{align*}
    (iii) It holds that
    \begin{align*}
        \frac{c_1^0}{c_2^0}
        =&l(m^*)^2\frac{(2kl+1)!_{(l)}}{(2k+1)!}
         \frac{ \frac{2}{l^{ k+1+\frac{1}{2} }} \frac{\Gamma\left(k+1+\frac{1}{2}\right)
        \Gamma\left(2k+1+\frac{3}{2l}\right)}
        {\Gamma\left(3k+2+\frac{1}{2}+\frac{3}{2l}\right)}
        }{
        \frac{2}{l^{ k+1+\frac{1}{2} }} 
        \frac{\Gamma\left( k+1+\frac{1}{2} \right)
        \Gamma\left(\frac{3}{2l}\right)}
        {\Gamma\left(k+1+\frac{1}{2}+\frac{3}{2l}\right)}
        }
    \\
       =&l\left(\frac{(2k+1)!!}{(2kl+1)!_{(2l)}}\right)^2
       \frac{(2kl+1)!_{(l)}}{(2k+1)!}
       \frac{\left(2k+\frac{3}{2l}\right)!}
       {\left(3k+\frac{3}{2}+\frac{3}{2l}\right)
       \left(3k+\frac{1}{2}+\frac{3}{2l}\right)\cdots
       \left(k+\frac{3}{2}+\frac{3}{2l}\right)}
       \\
       =&\left( \prod_{j=1}^{k}
       \frac{(2j+1)^2(2jl+1)((2j-1)l+1)}{(2jl+1)^2(2j+1)(2j)}
       \right)
       \cdot \left(
       l\prod_{j=0}^{2k}
       \frac{j+\frac{3}{2l}}{j+k+\frac{3}{2}+\frac{3}{2l}}
       \right)
       \\
       =&\left( \prod_{j=1}^{k}
       \left(1-\frac{l-1}{2j(2jl+1)}\right)
       \right)
       \cdot \left(
       l\prod_{j=0}^{2k}
       \frac{j+\frac{3}{2l}}{j+k+\frac{3}{2}+\frac{3}{2l}}
       \right)
       \\
       <&l\prod_{j=0}^{2k}
       \frac{j+\frac{3}{2l}}{j+k+\frac{3}{2}+\frac{3}{2l}}
       \\
       =&\frac{\frac{3}{2}}{k+\frac{3}{2}+\frac{3}{2l}}
       \prod_{j=1}^{2k}
       \frac{j+\frac{3}{2l}}{j+k+\frac{3}{2}+\frac{3}{2l}}
       \\
       =&\frac{\frac{3}{2}}{k+\frac{3}{2}+\frac{3}{2l}}
       \prod_{j=1}^{2k} \left(
       \frac{j-2k+\frac{3}{2l}-3}
       {3\left(j+k+\frac{3}{2}+\frac{3}{2l}\right)}+\frac{2}{3}
       \right)
       <\frac{3}{5} \left(\frac{2}{3}\right)^{2k},
    \end{align*}
    and when $1\le i\le 2k+1$
    \begin{align*}
        \frac{b_2^i}{c_2^0}=&m^*\frac{((i-1)l+1)!_{(l)}}{i!}
        \frac{ \frac{2}{l^{i+\frac{1}{2} }} \frac{\Gamma\left(i+\frac{1}{2}\right)
        \Gamma\left(k+1+\frac{3}{2l}\right)}
        {\Gamma\left(k+i+1+\frac{1}{2}+\frac{3}{2l}\right)}
        }{
        \frac{2}{l^{ k+1+\frac{1}{2} }} 
        \frac{\Gamma\left( k+1+\frac{1}{2} \right)
        \Gamma\left(\frac{3}{2l}\right)}
        {\Gamma\left(k+1+\frac{1}{2}+\frac{3}{2l}\right)}
        }
        \\
        =&l^{k+1-i}\frac{(2k+1)!!}{(2kl+1)!_{(2l)}}
        \frac{((i-1)l+1)!_{(l)}}{i!}
        \\
        &\times
        \frac{\left(k+\frac{3}{2l}\right)!}
        {\left(k+\frac{1}{2}\right)\cdots\left(i+\frac{1}{2}\right) }
        \frac{1}{ \left(k+i+\frac{1}{2}+\frac{3}{2l}\right) \cdots
        \left(k+\frac{3}{2}+\frac{3}{2l}\right) }
        \\
        =&\frac{(2kl+3)!_{(2l)}}{(2kl+1)!_{(2l)}}
        \frac{((i-1)l+1)!_{(l)}}{(il)!_{(l)}}
        \frac{\left(i-\frac{1}{2}\right)\cdots\frac{1}{2}}
        {\left(k+i+\frac{1}{2}+\frac{3}{2l}\right) \cdots
        \left(k+\frac{3}{2}+\frac{3}{2l}\right)}
        \\
        \le&\frac{1}{l}\frac{(2kl+3)!_{(2l)}}{(2kl+1)!_{(2l)}}
        \frac{\left(i-\frac{1}{2}\right)\cdots\frac{1}{2}}
        {\left(k+i+\frac{1}{2}+\frac{3}{2l}\right) \cdots
        \left(k+\frac{3}{2}+\frac{3}{2l}\right)}
        \\
        =&\frac{1}{l}\prod_{j=0}^{k}\frac{2jl+3}{2jl+1}
        \cdot\left(\frac{1}{2k+3+\frac{3}{l}}
        \prod_{j=2}^{i}\frac{j-\frac{1}{2}}{k+j+\frac{1}{2}+\frac{3}{2l}}
        \right)
        \\
        =&\frac{1}{l}\prod_{j=0}^{k}\frac{2jl+3}{2jl+1}
        \cdot\left(\frac{1}{2k+3+\frac{3}{l}}
        \prod_{j=2}^{i}\left(
        \frac{j-2k-\frac{5}{2}-\frac{3}{l}}
        {3\left(k+j+\frac{1}{2}+\frac{3}{2l}\right)} + \frac{2}{3}
        \right) \right)
        \\
        \le&\frac{1}{l}\prod_{j=0}^{k}\frac{2jl+3}{2jl+1}
        \cdot\left(\frac{1}{2k+3+\frac{3}{l}}
        \prod_{j=2}^{i}\left(
        \frac{-\frac{3}{2}-\frac{3}{l}}
        {3\left(k+j+\frac{1}{2}+\frac{3}{2l}\right)} + \frac{2}{3}
        \right) \right)
        \\
        \le&\frac{1}{(2k+3)l+3}\prod_{j=0}^{k}\frac{2jl+3}{2jl+1}
        \cdot
        \left( \frac{2}{3} \right) ^{i-1}
        \\
        =&\frac{1}{l+1}\prod_{j=1}^{k}
        \frac{((2j+1)l+3)(2jl+3)}{((2j+3)l+3)(2jl+1)}
        \cdot
        \left( \frac{2}{3} \right) ^{i-1}
        \\
        =&\frac{1}{l+1} 
        \prod_{j=1}^{k} \left( 1-
        \frac{4jl(l-1)-6}{(2jl+3l+3)(2jl+1)}
        \right)
        \cdot \left( \frac{2}{3} \right) ^{i-1}
        \\
        \le&\frac{1}{l+1} \frac{(3l+3)(2l+3)}{(5l+3)(2l+1)}
        \left( \frac{2}{3} \right) ^{i-1}
        \\
         =&\frac{3}{5l+3}
         \left(1+\frac{2}{2l+1}\right)
        \left( \frac{2}{3} \right) ^{i-1}
        \le\frac{21}{65}\left( \frac{2}{3} \right) ^{i-1}.
    \end{align*}
    Then,
    \begin{align*}
        -c_1^0-\sum_{i=1}^{2k+1}b_2^i +c_2^0=&
        c_2^0\left(1-\frac{c_1^0}{c_2^0}
        -\sum_{i=1}^{2k+1}\frac{b_2^i}{c_2^0}\right)
        \\
        >&c_2^0\left(1- 
        \frac{3}{5} \left(\frac{2}{3}\right)^{2k}
        -\frac{21}{65}\sum_{i=1}^{2k+1}\left( \frac{2}{3} \right) ^{i-1}
        \right)
        \\
        =&c_2^0 \left(
        \frac{2}{65}+\frac{9}{130}
        \left(\frac{2}{3}\right)^{2k+1}
        \right)>0.
    \end{align*}
    
\end{proof}
\end{prop}

\begin{lem}\label{Lemma: estimate of main terms for l>=3}
    $b_1^0-b_2^0-b_4^0>0,$
    when $l\ge3.$

    \begin{proof}
    It holds that
\begin{align}
    \frac{b_1^0}{b_2^0}=&
    l(m^*)^2\frac{(2kl+1)!_{(l)}}{(2k+1)!}
    \frac{ \frac{2}{l^{ \frac{1}{2} }} \frac{\Gamma\left(\frac{1}{2}\right)
        \Gamma\left(3k+2+\frac{3}{2l}\right)}
        {\Gamma\left(3k+2+\frac{1}{2}+\frac{3}{2l}\right)}
        }{
        \frac{2}{l^{ \frac{1}{2} }} 
        \frac{\Gamma\left( \frac{1}{2} \right)
        \Gamma\left(k+1+\frac{3}{2l}\right)}
        {\Gamma\left(k+1+\frac{1}{2}+\frac{3}{2l}\right)}
        }
    \notag    \\
        =&l\left(\frac{(2k+1)!!}{(2kl+1)!_{(2l)}}\right)^2\frac{(2kl+1)!_{(l)}}{(2k+1)!}
        \frac{\left(3k+1+\frac{3}{2l}\right)\left(3k+\frac{3}{2l}\right)
        \cdots\left(k+1+\frac{3}{2l}\right)}
        {\left(3k+\frac{3}{2}+\frac{3}{2l}\right)\left(3k+\frac{1}{2}+\frac{3}{2l}\right)
        \cdots\left(k+\frac{3}{2}+\frac{3}{2l}\right)}
     \notag   \\
        =&l\prod_{j=1}^{k}
        \left(\frac{2j+1}{2jl+1}\right)^2\frac{(2jl+1)((2j-1)l+1)}{(2j+1)(2j)}
        \cdot \prod_{j=1}^{2k+1}
        \frac{k+j+\frac{3}{2l}}{k+j+\frac{1}{2}+\frac{3}{2l}}
    \notag    \\
        =&l\prod_{j=1}^{k}
        \frac{2j+1}{2j}\frac{2j-1+\frac{1}{l}}{2j+\frac{1}{l}}
        \cdot \prod_{j=1}^{2k+1}
        \frac{k+j+\frac{3}{2l}}{k+j+\frac{1}{2}+\frac{3}{2l}}
    \label{b_1^0/b_2^0 (a)} \\
    >&l\prod_{j=1}^{k}
        \frac{2j+1}{2j}\frac{2j-1}{2j}
        \cdot \prod_{j=1}^{2k+1}
        \frac{k+j}{k+j+\frac{1}{2}}
    \notag \\
    =&l\prod_{j=1}^{k}
        \frac{4j^2-1}{4j^2}
        \cdot \frac{2}{3}\prod_{j=1}^{k}
        \frac{j+\frac{1}{2}}{j}
        \frac{(3j+1)(3j)(3j-1)}
        {\left(3j+\frac{3}{2}\right)\left(3j+\frac{1}{2}\right)
        \left(3j-\frac{1}{2}\right)}
    \notag \\
    =&\frac{2l}{3}\prod_{j=1}^{k}
        \frac{4j^2-1}{4j^2}
        \cdot \prod_{j=1}^{k}
        \frac{9j^2-1}{9j^2-\frac{1}{4}}
    \label{b_1^0/b_2^0 (b)} .
\end{align}
On the other hand, when $l\ge3,$ we have
\begin{align}
    \frac{b_4^0}{b_2^0}=&
    lm^*\frac{(k+1)!((k-1)l+1)!_{(l)}}{(2k+1)!}
    \frac{ \frac{2}{l^{ \frac{1}{2} }} \frac{\Gamma\left(\frac{1}{2}\right)
        \Gamma\left(2k+1+\frac{3}{2l}\right)}
        {\Gamma\left(2k+1+\frac{1}{2}+\frac{3}{2l}\right)}
        }{
        \frac{2}{l^{ \frac{1}{2} }} 
        \frac{\Gamma\left( \frac{1}{2} \right)
        \Gamma\left(k+1+\frac{3}{2l}\right)}
        {\Gamma\left(k+1+\frac{1}{2}+\frac{3}{2l}\right)}
        }
    \notag \\
    =&lm^*\frac{(k+1)!((k-1)l+1)!_{(l)}}{(2k+1)!}
    \frac{\left(2k+\frac{3}{2l}\right)\left(2k-1+\frac{3}{2l}\right)
    \cdots\left(k+1+\frac{3}{2l}\right)}
    {\left(2k+\frac{1}{2}+\frac{3}{2l}\right)
    \left(2k-\frac{1}{2}+\frac{3}{2l}\right)
    \cdots\left(k+\frac{3}{2}+\frac{3}{2l}\right)}
    \label{b_4^0/b_2^0 (a)} \\
    \le&lm^*\frac{(k+1)!((k-1)l+1)!_{(l)}}{(2k+1)!}
    \frac{\left(2k+\frac{1}{2}\right)\left(2k-\frac{1}{2}\right)
    \cdots\left(k+\frac{3}{2}\right)}
    {\left(2k+1\right)
    \left(2k\right)
    \cdots\left(k+2\right)}
    \notag \\
    =&l\prod_{j=1}^k\frac{(j+1)(jl-l+1)}{(2jl+1)(2j)}
    \cdot \prod_{j=1}^k
    \frac{(j+1)\left(2j+\frac{1}{2}\right)\left(2j-\frac{1}{2}\right)}
    {\left(j+\frac{1}{2}\right)(2j+1)(2j)}
    \notag \\
    =&l\prod_{j=1}^k\frac{(j+1)(jl-l+1)}{(2jl+1)(2j)}
    \cdot \prod_{j=1}^k \left( 1-
    \frac{5j+1}{4j(2j+1)^2} \right)
    \notag \\
    \le&\frac{5l}{6}\prod_{j=1}^k\frac{(j+1)(jl-l+1)}{(2jl+1)(2j)}
    \notag \\
    =&\frac{5}{3\cdot2^k}\frac{l}{2l+1}\prod_{j=2}^k
    \frac{j^2l+j-l+1}{2j^2l+j}
    \le\frac{5}{3\cdot2^k}\frac{l}{2l+1}
    \le\frac{5l}{21\cdot2^k}. 
    \label{b_4^0/b_2^0 (b)}
\end{align}
Now for $k=1,2,$  by \eqref{b_1^0/b_2^0 (b)} and \eqref{b_4^0/b_2^0 (b)}, 
we obtain 
\begin{align*}
    \left.\left(\frac{b_1^0}{b_2^0}-\frac{b_4^0}{b_2^0}\right)\right|_{k=1}
    &\ge l\left(\frac{16}{35}-\frac{5}{42}\right)
    \ge 3\cdot\frac{71}{210}>1
    \\
    \left.\left(\frac{b_1^0}{b_2^0}-\frac{b_4^0}{b_2^0}\right)\right|_{k=2}
    &\ge l\left(\frac{60}{143}-\frac{5}{84}\right)
    \ge 3\cdot\frac{4325}{12012}>1.
\end{align*}
Further, by Lemma \ref{Lemma: infinite product converge}
with some simple calculation, Eq. \eqref{b_1^0/b_2^0 (b)} becomes
\begin{align*}
    \frac{b_1^0}{b_2^0}
    >&\frac{2l}{3}\prod_{j=1}^{+\infty} \left( 
        1 - \frac{1}{4j^2} \right)
        \cdot \prod_{j=1}^{+\infty} \left(
        1- \frac{3}{36j^2-1} \right)
    \notag \\ \notag
    >&\frac{2l}{3}
    \exp\left(-\frac{1}{3N}\right)
    \prod_{j=1}^{N} \left( 
        1 - \frac{1}{4j^2} \right)
    \cdot \exp\left(-\frac{3}{32N}\right)
    \prod_{j=1}^{N} \left(
        1- \frac{3}{36j^2-1} \right)
\end{align*}
for any $N\ge1.$  
Taking $N=10,$ we find that for $l,k\ge3,$
\begin{align}
    \frac{b_1^0}{b_2^0}-\frac{b_4^0}{b_2^0}
    >&l\left(\frac{2}{3}\exp{\left(-\frac{41}{960}\right)}
    \prod_{j=1}^{10}\left(1-\frac{1}{4j^2}\right)
    \left(1-\frac{3}{36j^2-1}\right)-\frac{5}{21\cdot2^k}
    \right)
    \notag\\
    \ge&3\left(\frac{2}{3}\exp{\left(-\frac{41}{960}\right)}
    \prod_{j=1}^{10}\left(1-\frac{1}{4j^2}\right)
    \left(1-\frac{3}{36j^2-1}\right)-\frac{5}{21\cdot2^3}
    \right)
    \notag \\
    \approx&3\cdot0.3338>1.
    \label{approximate (1)}
\end{align}
    \end{proof}
\end{lem}

\begin{prop}\label{Lemma: estimate of main terms}
    $b_1^0-b_2^0+\sum_{i=1}^{2k+1}c_2^i-b_4^0>0.$
    \begin{proof}
    By Lemma \ref{Lemma: estimate of main terms for l>=3}, 
    only the case $l=2$ is considered.
    When $l=2,$ by Lemma \ref{Lemma: infinite product converge}
    with some simple calculation, Eq. \eqref{b_1^0/b_2^0 (a)} becomes
    \begin{align}
        \frac{b_1^0}{b_2^0}=&2\prod_{j=1}^{k}
        \frac{2j+1}{2j}\frac{2j-1+\frac{1}{2}}{2j+\frac{1}{2}}
        \cdot \prod_{j=1}^{2k+1}
        \frac{k+j+\frac{3}{4}}{k+j+\frac{1}{2}+\frac{3}{4}}
        \notag\\
        =&\frac{14}{9}
        \prod_{j=1}^k \left(
        1-\tfrac{1}{2j(4j+1)}
        \right)
        \cdot
        \prod_{j=1}^k \left(
        1-\tfrac{48(48j^2+36j+5)}{(4j+3)(12j+9)(12j+5)(12j+1)}
        \right)
        \notag\\
        >&\frac{14}{9} 
        \prod_{j=1}^{+\infty} \left(
        1-\tfrac{1}{2j(4j+1)}
        \right)
        \notag\cdot
        \prod_{j=1}^{+\infty} \left(
        1-\tfrac{48(48j^2+36j+5)}{(4j+3)(12j+9)(12j+5)(12j+1)}
        \right)
        \notag\\
        \begin{split} \label{b_1^0/b_2^0 (c)}
        >&\frac{14}{9} \exp{\left(-\tfrac{1}{8N}\right)}
        \prod_{j=1}^{N} \left(
        1-\tfrac{1}{2j(4j+1)}
        \right)
        \\
        &\times\exp{\left(-\tfrac{1}{3N}\right)}
        \prod_{j=1}^{N} \left(
        1-\tfrac{48(48j^2+36j+5)}{(4j+3)(12j+9)(12j+5)(12j+1)}
        \right)
        \end{split}
        \\ 
        \triangleq& \frac{14}{9} \lambda_N
        \notag
    \end{align}
    for any $N\ge1,$
    and Eq. \eqref{b_4^0/b_2^0 (a)} becomes
    \begin{align}
        \frac{b_4^0}{b_2^0}=&2m^*\frac{(k+1)!(2k-1)!!}{(2k+1)!}
        \frac{\left(2k+\frac{3}{4}\right)\left(2k-\frac{1}{4}\right)
        \cdots\left(k+\frac{7}{4}\right)}
        {\left(2k+\frac{5}{4}\right)
        \left(2k+\frac{1}{4}\right)
        \cdots\left(k+\frac{9}{4}\right)}
        \notag \\
        =&2\prod_{j=1}^k\frac{(j+1)(2j-1)}{(4j+1)(2j)}
        \cdot \prod_{j=1}^k 
        \frac{\left(j+\frac{5}{4}\right)\left(2j+\frac{3}{4}\right)
        \left(2j-\frac{1}{4}\right)}
        {\left(j+\frac{3}{4}\right)\left(2j+\frac{5}{4}\right)
        \left(2j+\frac{1}{4}\right)}
        \notag \\
        =&\frac{1}{2^{k-1}} \prod_{j=1}^k 
        \frac{2j^2+j-1}{4j^2+j}\
        \cdot
        \prod_{j=1}^k \left( 1-
        \frac{6(16j+5)}{(4j+3)(8j+5)(8j+1)} \right)
        \notag \\
        \le& \frac{1}{2^{k-1}} 
        \cdot 
        \frac{2}{5}
        \cdot
        \frac{11}{13}
        =\frac{11}{65\cdot 2^{k-2}}.
        \label{b_4^0/b_2^0 (c)}
    \end{align}
    On the other hand, 
        \begin{align*}
            \frac{c_2^i}{b_2^0}=&\frac{1}{m^*}\frac{((i-1)l+1)!_{(l)}}{i!}
            \frac{ \frac{2}{l^{ k+i+1+\frac{1}{2} }} \frac{\Gamma\left(k+i+1+\frac{1}{2}\right)
            \Gamma\left(\frac{3}{2l}\right)}
            {\Gamma\left(k+i+1+\frac{1}{2}+\frac{3}{2l}\right)}
            }{
            \frac{2}{l^{ \frac{1}{2} }} 
            \frac{\Gamma\left( \frac{1}{2} \right)
            \Gamma\left(k+1+\frac{3}{2l}\right)}
            {\Gamma\left(k+1+\frac{1}{2}+\frac{3}{2l}\right)}
            }
        \\
            =&\frac{(2kl+1)!_{(2l)}}{(2kl+l)!_{(2l)}}
            \frac{((i-1)l+1)!_{(l)}}{(il)!_{(l)}}
            \frac{\left(k+i+\frac{1}{2}\right)!}
            {\left(k+\frac{3}{2l}\right)!}
        \\
            &\times\frac{1}
            {\left(k+i+\frac{1}{2}+\frac{3}{2l}\right)\left(k+i-\frac{1}{2}+\frac{3}{2l}\right)
            \cdots\left(k+\frac{3}{2}+\frac{3}{2l}\right)}
        \\
            =&\frac{(2kl+1)!_{(2l)}}{(2kl+3)!_{(2l)}}
            \frac{((i-1)l+1)!_{(l)}}{(il)!_{(l)}}
            \prod_{j=k+1}^{k+i} \frac{j+\frac{1}{2}}{j+\frac{1}{2}+\frac{3}{2l}}.
        \end{align*}
        When $l=2,$ it becomes
        $$
          \frac{c_2^i}{b_2^0}= 
          \prod_{j=0}^k \frac{4j+1}{4j+3}
          \cdot 
          \prod_{j=1}^{i} \frac{2j-1}{2j}
          \cdot
          \prod_{j=k+1}^{k+i} \frac{4j+2}{4j+5}.
        $$
        We then estimate the three rational product terms above as follows:
        \begin{align*}
            (k+1)^{\frac{1}{2}} \prod_{j=0}^k \frac{4j+1}{4j+3}
            =&\frac{1}{3} \sqrt{ \prod_{j=1}^k
            \frac{j+1}{j}\left(\frac{4j+1}{4j+3}\right)^2 }
            \\
            =&\frac{1}{3} \sqrt{ \prod_{j=1}^k
            \left(1+\frac{1}{j(4j+3)^2}\right) }
            \ge \frac{5}{21}\sqrt{2},
        \end{align*} 
        \begin{align*}
            (i+1)^{\frac{1}{2}}\prod_{j=1}^{i} \frac{2j-1}{2j}
            =& \sqrt{\prod_{j=1}^i
            \left( 1-\frac{3j-1}{4j^3} \right) }
            \\
            >&\sqrt{\prod_{j=1}^{+\infty}
            \left( 1-\frac{3j-1}{4j^3} \right) }
            \\
            >&\exp{\left(-\tfrac{1}{2N}\right)} 
            \sqrt{\prod_{j=1}^{N}
            \left( 1-\tfrac{3j-1}{4j^3} \right) }
            \triangleq \mu_N
        \end{align*}
        for any $N\ge1,$ and
        \begin{align*}
            \left(\tfrac{k+i+1}{k+1}\right)^{\frac{3}{4}} 
            \prod_{j=k+1}^{k+i} \tfrac{4j+2}{4j+5}
            =&\sqrt[4]{\prod_{j=k+1}^{k+i}
            \left( 1+
            \tfrac{288j^5+1072j^4+1439j^3+816j^2+176j+16}
            {j^3(4j+5)^4}
            \right)}
            \\
            >&1.
        \end{align*}
        Now, it follows that
        \begin{align}
            \sum_{i=1}^{2k+1} \frac{c_2^i}{b_2^0} 
            >& \sum_{i=1}^{2k+1} \left(
            \tfrac{5\sqrt{2}}{21} (k+1)^{-\frac{1}{2}}
            \cdot
            \mu_N (i+1)^{-\frac{1}{2}}
            \cdot
            \left(\tfrac{k+i+1}{k+1}\right)^{-\frac{3}{4}} 
            \right)
            \notag \\
            >&\frac{5\sqrt{2}}{21} \mu_N 
            (k+1)^{\frac{1}{4}} \sum_{i=1}^{2k+1}
            (i+1)^{-\frac{1}{2}} (k+i+2)^{-\frac{3}{4}}
            \notag \\
            =&\frac{5\sqrt{2}}{21} \mu_N \left(
            \frac{1}{k+1} \sum_{i=1}^{2k+1}
            \left(\frac{i+1}{k+1}\right)^{-\frac{1}{2}}
            \left( 1 + \frac{i+1}{k+1} \right)^{-\frac{3}{4}}
            \right)
            \notag \\
            >&\frac{5\sqrt{2}}{21} \mu_N
            \int_{\frac{2}{k+1}}^{2+\frac{1}{k+1}}
            x^{-\frac{1}{2}}(1+x)^{-\frac{3}{4}}
            \mathrm{d}x
            \notag \\
            >&\frac{5\sqrt{2}}{21} \mu_N
            \int_{\frac{2}{k+1}}^{2}
            x^{-\frac{1}{2}}(1+x)^{-\frac{3}{4}}
            \mathrm{d}x.
            \label{sum c_2^i/b_2^0}
        \end{align}
        In terms of \eqref{b_1^0/b_2^0 (c)}, \eqref{b_4^0/b_2^0 (c)}
        and \eqref{sum c_2^i/b_2^0}, we obtain
        \begin{align*}
            \tfrac{b_1^0}{b_2^0}-1
            +\sum_{i=1}^{2k+1} \tfrac{c_2^i}{b_2^0}
            -\tfrac{b_4^0}{b_2^0}
            >& \tfrac{14}{9}\lambda_N
            +\tfrac{5\sqrt{2}}{21} \mu_N
            \int_{\frac{2}{k+1}}^{2}
            x^{-\frac{1}{2}}(1+x)^{-\frac{3}{4}}
            \mathrm{d}x
            -\tfrac{11}{65\cdot 2^{k-2}}-1,
        \end{align*}
        the right-hand side of which is increasing with respect to $k$ for any fixed $N.$
        If take $N=3,$ we find that for $k\ge3$
        \begin{align}
            \tfrac{b_1^0}{b_2^0}-1
            +\sum_{i=1}^{2k+1} \tfrac{c_2^i}{b_2^0}
            -\tfrac{b_4^0}{b_2^0}
            >& \tfrac{14}{9}\lambda_3
            +\tfrac{5\sqrt{2}}{21} \mu_3
            \int_{\frac{1}{2}}^{2}
            x^{-\frac{1}{2}}(1+x)^{-\frac{3}{4}}
            \mathrm{d}x
            -\tfrac{11}{130}-1
            \notag\\
            \approx& 0.0147>0.
            \label{approximate (2)}
        \end{align}
        When $k=1,2,$ by direct calculation, we see that
        \begin{align*}
            \left.\left(
            \tfrac{b_1^0}{b_2^0}-1
            +\sum_{i=1}^{2k+1} \tfrac{c_2^i}{b_2^0}
            -\tfrac{b_4^0}{b_2^0}
            \right) \right|_{l=2,k=1} =& \frac{1531}{23205},
            \\
            \left.\left(
            \tfrac{b_1^0}{b_2^0}-1
            +\sum_{i=1}^{2k+1} \tfrac{c_2^i}{b_2^0}
            -\tfrac{b_4^0}{b_2^0}
            \right) \right|_{l=2,k=2} =& \frac{203341}{759220}.
        \end{align*}
        Thus, it is indeed positive for all $k,$ 
        which completes the proof.
    \end{proof}
\end{prop}

\begin{rmk}
    \normalfont{
    In the proofs of Lemma \ref{Lemma: estimate of main terms for l>=3}
    and Proposition \ref{Lemma: estimate of main terms},
    precisely \eqref{approximate (1)} and \eqref{approximate (2)},
    we use numerical approximations of exponent functions and an integral,
    i.e., $
           \int_{\frac{1}{2}}^{2}
            x^{-\frac{1}{2}}(1+x)^{-\frac{3}{4}}
            \mathrm{d}x,
          $
    though, they can be estimated to formulas 
    which are the sum of $n$th-roots of rational numbers.
    For example, for $x>0$ and sufficiently large $M$
    \begin{align*}
        \exp(-x)=\sum_{i=0}^{\infty}(-1)^{i}\frac{x^i}{i!}
        >\sum_{i=0}^{M} \left(
        \frac{x^{2i}}{(2i)!}-\frac{x^{2i+1}}{(2i+1)!}
        \right),
    \end{align*}
    and
    \begin{align*}
        \int_{\frac{1}{2}}^{2}
            x^{-\frac{1}{2}}(1+x)^{-\frac{3}{4}}
            \mathrm{d}x
        >\frac{1}{M} \sum_{i=1}^{M}
            \left(\frac{3i}{2M}+\frac{1}{2}\right)^{-\frac{1}{2}}
            \left( \frac{3i}{2M}+\frac{3}{2} \right)^{-\frac{3}{4}}.
    \end{align*}
    Therefore, taking $M$ large enough,
    the numerical approximations in
    \eqref{approximate (1)} and \eqref{approximate (2)}
    can be replaced by a rational number and
    the sum of $n$th-roots of rational numbers, respectively,
    which makes our proofs more theoretical.
    }
\end{rmk}

Now, by Propositions \ref{Lemma: estimate of infinitesimal} and
\ref{Lemma: estimate of main terms}, we conclude that
for any $l(\ge2),k\in\mathbb{N}^*$
\begin{align*}
    V=&b_1^0-c_1^0
     -\sum_{i=0}^{2k+1}\left(b_{2}^i-c_{2}^i\right)
    + \sum_{i=0}^{k}\left(b_{3}^i-c_{3}^i\right)
    - \sum_{i=0}^{k+1}\left(b_{4}^i-c_{4}^i\right)
    \\
    =&\left(b_1^0-b_2^0+\sum_{i=1}^{2k+1}c_2^i-b_4^0\right)
    +\left(-c_1^0-\sum_{i=1}^{2k+1}b_2^i+c_2^0\right)
    \\
    &+\sum_{i=0}^{k}\left(b_3^i-b_4^{i+1}\right)
    +\sum_{i=0}^{k}\left(-c_3^i+c_4^{i+1}\right)
    +c_4^{0}>0.
\end{align*}

\subsection{Proof of $W>0$}

Similarly, $W$ is divided into six parts:
$W=W_1-W_2+W_3-W_4+W_5-W_6,$
where
\begin{align*}
  &W_1=l^2\left(m^*\right)^3
  \int_0^\Omega\text{Sn}^2(\theta)\text{Cs}^{2(3k+1)l+2}(\theta)
  \mathrm{d}\theta,
  \\
  &W_2=\int_0^\Omega
  \text{Sn}^{6k+4}(\theta)\text{Cs}^2(\theta)
  \mathrm{d}\theta,
  \\
  &W_3=2l(m^*)^2 \int_0^\Omega
  \text{Sn}^{2k+2}(\theta)\text{Cs}^{2(2k+1)l+2}(\theta)
  \mathrm{d}\theta,
  \\
  &W_4=(lm^*)^2\int_0^\Omega
  \text{Sn}^{2k+4}(\theta)\text{Cs}^{4kl+2}(\theta)
  \mathrm{d}\theta,
  \\
  &W_5=m^* \int_0^\Omega
  \text{Sn}^{4k+2}(\theta)\text{Cs}^{2(k+1)l+2}(\theta)
  \mathrm{d}\theta,
  \\
  &W_6=2lm^* \int_0^\Omega
  \text{Sn}^{4k+4}(\theta)\text{Cs}^{2kl+2}(\theta)
  \mathrm{d}\theta.
\end{align*}
We are going to consider
$W_1-W_2$
and
$W_3-W_4+W_5-W_6,$
respectively.

Again, with the help of Lemma \ref{Lemma: infinite product converge},
we have
\begin{align}
  \frac{W_2}{W_1} =&\frac{1}{l^2(m^*)^3}
            \frac{ \frac{2}{l^{ 3k+2+\frac{1}{2} }} 
            \frac{ \Gamma\left(3k+2+\frac{1}{2}\right)\Gamma\left(\frac{3}{2l}\right) }
            {\Gamma\left(3k+2+\frac{1}{2}+\frac{3}{2l}\right)}
            }{
            \frac{2}{l^{ \frac{3}{2} }} 
            \frac{ \Gamma\left( \frac{3}{2} \right)\Gamma\left(3k+1+\frac{3}{2l}\right) }
            {\Gamma\left(3k+2+\frac{1}{2}+\frac{3}{2l}\right)}
            }
  \notag \\
  =&\frac{1}{l^{3k+3}} \left(\frac{(2kl+1)!_{(2l)}}{(2k+1)!!}\right)^3
  \frac{\left(3k+\frac{3}{2}\right)\left(3k-1+\frac{3}{2}\right)\cdots\frac{3}{2}}
  {\left(3k+\frac{3}{2l}\right)\left(3k-1+\frac{3}{2l}\right)\cdots\frac{3}{2l}}
  \notag \\
  =&\frac{1}{l^2} \prod_{j=1}^k
  \frac{\left(j+\frac{1}{2l}\right)^3
  \left(3j+\frac{3}{2}\right)\left(3j+\frac{1}{2}\right)\left(3j-\frac{1}{2}\right) }
  {\left(j+\frac{1}{2}\right)^3
  \left(3j+\frac{3}{2l}\right)\left(3j-1+\frac{3}{2l}\right)\left(3j-2+\frac{3}{2l}\right) }
  \notag \\
  =&\frac{1}{l^2} \prod_{j=1}^k \left(
  \frac{j+\frac{1}{2l}}{3\left(j-\frac{1}{3}+\frac{1}{2l}\right)}
  \cdot
  \frac{j+\frac{1}{2l}}{3\left(j-\frac{2}{3}+\frac{1}{2l}\right)}
  \cdot
  \frac{ \left(3j+\frac{1}{2}\right)\left(3j-\frac{1}{2}\right) }
  { \left(j+\frac{1}{2}\right)^2 }
  \right)
  \notag \\
  <&\frac{1}{4} \prod_{j=1}^k \left(
  \frac{j}{3\left(j-\frac{1}{3}\right)}
  \cdot
  \frac{j}{3\left(j-\frac{2}{3}\right)}
  \cdot
  \frac{ \left(3j+\frac{1}{2}\right)\left(3j-\frac{1}{2}\right) }
  { \left(j+\frac{1}{2}\right)^2 }
  \right)
  \label{W2/W1 (a)}
  \\
  =&\frac{1}{4} \prod_{j=1}^{k} \left(
  1+\frac{18j^2+j-2}{(2j+1)^2(3j-1)(3j-2)}
  \right)
  \notag \\
  <&\frac{1}{4} \prod_{j=1}^{+\infty} \left(
  1+\frac{18j^2+j-2}{(2j+1)^2(3j-1)(3j-2)}
  \right)
  \notag \\
  <&\frac{1}{4}\exp{\left(\frac{1}{2}\right)} \prod_{j=1}^2
  \left(
  1+\frac{18j^2+j-2}{(2j+1)^2(3j-1)(3j-2)}
  \right)
  \notag \\
  =&\frac{1}{4}\cdot\frac{35}{18}\cdot\frac{143}{125}
  \cdot\exp{\left(\frac{1}{2}\right)}
  \approx0.9169<1.
  \label{W2/W1 (b)}
\end{align}
On the other hand, for the rest terms we have
\begin{align}
    &-\frac{W_4}{W_3}+\frac{W_5}{W_3}-\frac{W_6}{W_3}
    \notag\\
    =&-\frac{k+\frac{3}{2}}{4k+\frac{3}{l}}
    +\frac{1}{l^{k+1}} \left( \frac{k+\frac{3}{2l}}{4k+3} -1 \right)
    \frac{(2kl+1)!_{(2l)}}{(2k+1)!!}
    \frac{\left(2k+\frac{3}{2}\right)\left(2k+\frac{1}{2}\right)
    \cdots\left(k+\frac{3}{2}\right)}
    {\left(2k+\frac{3}{2l}\right)\left(2k-1+\frac{3}{2l}\right)
    \cdots\left(k+\frac{3}{2l}\right)}
    \notag\\
    =&-\frac{k+\frac{3}{2}}{4k+\frac{3}{l}}
    -\frac{3k+3-\frac{3}{2l}}{4k+3}
    \prod_{j=0}^{k}\frac{2j+\frac{1}{l}}{2j+1}
    \frac{k+j+\frac{3}{2}}{k+j+\frac{3}{2l}}
    \notag\\
    >&-\frac{2k+3}{8k}-\frac{3k+3}{4k+3}
    \prod_{j=0}^{k}\frac{2j+\frac{1}{l}}{2j+1}
    \frac{k+j+\frac{3}{2}}{k+j+\frac{3}{2l}}
    \notag\\
    \ge&-\frac{2k+3}{8k}-\frac{3k+3}{4k+3}
    \prod_{j=0}^{k}\frac{2j+\frac{1}{2}}{2j+1}
    \frac{k+j+\frac{3}{2}}{k+j+\frac{3}{4}}
    \notag\\
    =&-\frac{2k+3}{8k}-\frac{3k+3}{4k+3}
    \prod_{j=1}^k \left(
    1-\frac{64j^3+84j^2+20j}{(2j+1)^2(8j+3)(8j-1)}
    \right)\triangleq \nu_k,
    \label{-W4/W3+W5/W3-W6/W3}
\end{align}
where the second inequality is deduced from the monotonicity of  \\
$\prod_{j=0}^{k}(j+u)/(k+j+3u)$ 
with respect to 
$u=1/(2l).$
In fact, we have
\begin{align*}
    \prod_{j=0}^{k}\frac{j+u}{k+j+3u}=&
    \sqrt{ \prod_{j=0}^{k}\frac{j+u}{k+j+3u}
    \cdot \prod_{j=0}^{k}\frac{k-j+u}{2k-j+3u} }
    \\
    =&\sqrt{ \prod_{j=0}^{k}\frac{j(k-j)+ku+u^2}{(k+j)(2k-j)+9ku+9u^2} }
    \\
    =&\frac{1}{3^{k+1}} \sqrt{ \prod_{j=0}^{k} \left(
    1- \frac{2(k-2j)^2}{(k+j)(2k-j)+9ku+9u^2} \right) },
\end{align*}
which is increasing with respect to $u.$

One can see that
$\nu_k$
is increasing.
Hence, by \eqref{W2/W1 (b)} and \eqref{-W4/W3+W5/W3-W6/W3}
we have for $l,k\ge2$
\begin{equation*}
  W>W_1\left(1-\frac{W_2}{W_1}\right)
   +W_3\left(1+
   \nu_k
   \right)
   >W_3(1+\nu_2)
   =W_3\cdot\frac{333}{16720}
   >0.
\end{equation*}
Additionally, when 
$k=1,$
by using \eqref{W2/W1 (a)} instead of \eqref{W2/W1 (b)}
we also have for $l\ge2$
\begin{align*}
  \left.\frac{W}{W_3}\right|_{k=1}
  >&\left.\frac{W_1}{W_3}\right|_{k=1}
  \cdot
  \left( 1 - \left.\frac{W_2}{W_1}\right|_{k=1}\right)
  +1+\nu_1
  \\
  =&\frac{3l}{2} \cdot
  \left( 1 - \left.\frac{W_2}{W_1}\right|_{k=1}\right)
  -\frac{169}{616}
  \\
   >&3\cdot\left( 1 - \frac{35}{72} \right) -\frac{169}{616}
   =\frac{37}{24}-\frac{169}{616}>0.
\end{align*}
It follows that $W$ is always positive.

Consequently, we have verified
$u_{3K+1}(\Omega)>0$
for all $l(\ge2), k\in\mathbb{N}^*,$
and the proof of 
Theorem \ref{Theorem_main result: unstable focus}
is completed.

\section*{Acknownledgement}
This paper is supported by the National Natural Science
Foundation of China (No. 12171491).




\begin{thebibliography}{10}
\expandafter\ifx\csname url\endcsname\relax
  \def\url#1{\texttt{#1}}\fi
\expandafter\ifx\csname urlprefix\endcsname\relax\def\urlprefix{URL }\fi
\expandafter\ifx\csname href\endcsname\relax
  \def\href#1#2{#2} \def\path#1{#1}\fi

\bibitem{dumortier2006qualitative}
F.~Dumortier, J.~Llibre, J.~C. Art{\'e}s, Qualitative theory of planar differential systems, Vol.~2, Springer, 2006.

\bibitem{andreev1958investigation}
A.~F. Andreev, Investigation of the behaviour of the integral curves of a system of two differential equations in the neighbourhood of a singular point, Trans. Amer. Math. Soc 8 (1958) 187--207.

\bibitem{alvarez2005monodromy}
M.~{\'A}lvarez, A.~Gasull, Monodromy and stability for nilpotent critical points, Int. J. Bifurc. Chaos 15~(4) (2005) 1253--1265.

\bibitem{lyapunov1966stability}
A.~M. Lyapunov, Stability of Motion, Vol.~30 of Mathematics in Science and Engineering, Academic Press, New York-London, 1966.

\bibitem{gasull2001new}
A.~Gasull, J.~Torregrosa, A new algorithm for the computation of the lyapunov constants for some degenerated critical points, Nonlinear Anal. Theory 47 (2001) 4479--4490.

\bibitem{moussu1982symetrie}
R.~Moussu, Sym{\'e}trie et forme normale des centres et foyers d{\'e}g{\'e}n{\'e}r{\'e}s, Ergod. Theory Dyn. Syst. 2~(2) (1982) 241--251.

\bibitem{chen2022new}
H.~Chen, R.~Zhang, X.~Zhang, New criterions on stability and order of analytic nilpotent foci, J. Differ. Equ. 338 (2022) 352--371.

\bibitem{galeotti1989bifurcations}
M.~Galeotti, F.~Gori, Bifurcations and limit cycles in a family of planar polynomial dynamical systems, Rend. Semin. Mat. (Torino) 46 (1989) 31--58.

\bibitem{garcia2014bifurcation}
J.~D. Garc{\'\i}a-Salda{\~n}a, A.~Gasull, H.~Giacomini, Bifurcation diagram and stability for a one-parameter family of planar vector fields, J. Math. Anal. Appl. 413~(1) (2014) 321--342.

\bibitem{caubergh2015bifurcation}
M.~Caubergh, Bifurcation of the separatrix skeleton in some 1-parameter families of planar vector fields, J. Differ. Equ. 259~(3) (2015) 989--1013.

\end{thebibliography}


\bibliographystyle{elsarticle-num}

\end{document}